\documentclass{article}
\usepackage{amsfonts}
\usepackage{amsmath}
\begin{document}
\title{Symmetric representation of the elements of finite groups}
\author{Z. Hasan and \and A. Kasouha \\{\em Department of Mathematics, California
State University,} \\ {\em 5500 University Parkway,  San
Bernardino, California 92407}}

\maketitle

\begin{abstract}
Many finite groups, including all finite non-abelian simple
groups, can be symmetrically generated by involutions.  In this
paper we give an algorithm to symmetrically represent elements of
finite groups and to transform symmetrically represented elements
in terms of their permutation representation. In particular, we
represent elements of the group $U_3(3) : 2$, where $U_3(3)$ is
the projective special unitary group of $3 \times 3$ matrices with
entries in the field $\mathbb{F}_3$, as a permutation on fourteen
letters from the projective general linear group $PGL_2(7)$,
followed by a word of length at most two in the fourteen
involutory symmetric generators.
\end{abstract}
\section{Introduction}
\label{sec:intro} The two general methods for working with groups,
permutations and matrices, are inconvenient or unmanageable for
large finite groups and in particular for the larger sporadic
groups, for example for the Monster group, the least degree of a
matrix representation is 196881 and it takes $n^3$ operations to
multiply two matrices of dimension n and the least degree of any
permutation representation is $10^{20}$. Matrix multiplication for
matrices of large size is very time-consuming and, although MAGMA
and GAP handle permutations of reasonably large size quite
efficiently, recording and transmitting elements is inconvenient.
The main purpose of this paper is to give an alternative and more
efficient method for working with groups. Double coset enumeration
can be performed on groups that possess generating sets of
involutions. Curtis has constructed several sporadic groups by
manual double coset enumeration; for references see \cite{Cur3}
and \cite{Cur4}. Now any finite group generated by a conjugacy
class of involutions, and hence all finite non-abelian simple
groups have symmetric generating sets of involutions (see
\cite{Bray1}). It is this technique of double coset enumeration
that allows us to write elements in a much more concise manner. In
the first author's work with Curtis (see \cite{Cur4} ), every
element of J1, usually written as a permutation on 266 letters,
was written as an element of $L_2(11)$ followed by a word of
length at most 4 in the symmetric generators and a program was
written to manipulate elements of $J_1$ written in this short
form.
\section{Symmetric generation of a group}
Let $\cal G$ be a group and
\begin{center}
\mbox{$ {\cal T} = \{t_{1},t_{2},\ldots ,t_{n}\} \subseteq {\cal
G} $},
\end{center}
then we define
\begin{center}
$ \overline{{\cal T}} = \{ T_{1},T_{2},\ldots ,T_{n} \} $,
\end{center}
where \( T_{i} = \langle t_{i}  \rangle \), the cyclic subgroup
generated by $t_{i}$; we further define \mbox{${\cal N} = {\cal
N}_{\cal G}(\overline{{\cal T}})$}, the set normalizer in $\cal G$
of $\overline{{\cal T}}$.
\par
\noindent If the following two conditions hold:
\begin{enumerate}
\item[(i)] $ {\cal G} = \langle {\cal T} \rangle $, and
\item[(ii)] $\cal N$ permutes $\overline{{\cal T}}$ transitively, not
necessarily faithfully,
\end{enumerate}
then, following Curtis and Hasan \cite{Cur4}, we say that $\cal T$
is a {\em symmetric generating set}
 for $\cal G$ . In these circumstances we call $\cal N$ the
{\em control subgroup}. Note that (i) and (ii) imply that $\cal G$
is a homomorphic image of the (infinite) {\em progenitor}
\begin{center}
$m^{*n}$ : $\cal N$,
\end{center}
where $m^{\star n}$ represents a free product of $n$ copies of the
cyclic group $C_{m}$, $m$ being the order of $t_{i}$, and  $\cal
N$ is a group of automorphisms of $m^{\star n}$ which permutes the
$n$ cyclic subgroups by conjugation. Thus, for $\pi \in {\cal N}$,
we have
\begin{center}
$t_{i}^{\pi} = t_{j}^{r}$,
\end{center}
where $r$ is an integer coprime to $m$. Of course, if $m = 2$ then
$\cal N$ will simply act by conjugation as permutations of the $n$
involutory symmetric generators. Now, since by the above elements
of $\cal N$ can be gathered on the left, every element of the
progenitor can be represented as $\pi w$, where $\pi \in {\cal N}$
and $w$ is a word in the symmetric generators. Indeed this
representation is unique provided $w$ is simplified so that
adjacent symmetric generators are distinct. Thus any additional
relator by which we must factor the progenitor to obtain $\cal G$
must have the form
\begin{center}
$\pi w(t_{1},t_{2},..,t_{n})$,
\end{center}
where $\pi \in {\cal N}$ and $w$ is a word in {\cal T}. In the
next section we describe how a particular factor group
\begin{equation}
\frac{m^{*n}:N}{{\pi}_{1}w_{1},..,{\pi}_{s}w_{s}}
\end{equation}
may be identified.
\section{Manual double coset enumeration}
Since in this paper we are only concerned with involutory
symmetric generators we restrict our attention to the case $m=2$.
Thus we seek homomorphic images of the progenitor
\begin{center}
$2^{*n}:{\cal N}$,
\end{center}
(where $\cal N$ is now a transitive permutation group on $n$
letters), which act faithfully on {\cal N} and on the generators
of the free product. It is convenient to identify the $n$ free
generators and $\cal N$ with their respective images. Thus
\begin{equation}
\frac{2^{*n}:{\cal N}}{{\pi}_{1}w_{1},..,{\pi}_{s}w_{s}} \cong
\langle {\cal N}, {\cal T} \mid t_{i}^{2} = 1, t_{i}^{\pi} =
t_{\pi (i)}, \pi_{1}w_{1} = ..= \pi_{s}w_{s} = 1 \rangle ,
\end{equation}
where $\pi \in {\cal N}, t_{i} \in {\cal T}$. Following
\cite{Cur2} we are allowing $i$ to stand for the symmetric
generator $t_{i}$ in expressions such as the above relation. By a
slight abuse of notation we also allow $i$ to denote the coset
${\cal N}{t_{i}}$, \mbox{\em ij} the coset ${\cal
N}{t_{i}}{t_{j}}$ etc., and we write, for instance,
\begin{center}
$ij \sim k$ to mean ${\cal N}{t_{i}}{t_{j}}={\cal N}{t_{k}}$.
\end{center}
Writing $ij=k$ would be the much stronger statement that
${t_{i}}{t_{j}}={t_{k}}$. Now since
\begin{displaymath}
{t_{i}}{\pi}={\pi}{t_{{\pi}(i)}},
\end{displaymath}
(or $i{\pi}={\pi}{i^{\pi}}$ as we shall more commonly write), the
permutations involved in any element of $\cal G$ can be gathered
on the left. Thus any element of $\cal G$ can be written as a
permutation belonging to $\cal N$ followed by a word in the
symmetric generators. Indeed, as mentioned in the last section, in
the case of the progenitor itself this representation is unique
provided the obvious cancellations are performed. Thus, if
\mbox{${\cal N}g{\cal N}$} is a double coset of $\cal N$ in $\cal
G$, we have
\begin{displaymath}
{\cal N}g{\cal N}={\cal N}{\pi}w{\cal N}={\cal N}w{\cal N},
\end{displaymath}
where $g={\pi}w \in {\cal G}$, with ${\pi} \in {\cal N}$, and $w$
is a word in the $t_{i}$. We denote this double coset by
\mbox{$[w]$}, e.g. \mbox{$[01]$} denotes the double coset
\mbox{${\cal N}{t_{0}}{t_{1}}{\cal N}$}. The double coset ${\cal
N}e{\cal N}={\cal N}$, where $e$ is the identity element, is
denoted by [$\star$]. Furthermore we define
\begin{displaymath}
{{\cal N}^{i}}={{\cal C}_{{\cal N}}}(t_{i});\  {{\cal N}^{ij}}=
{{\cal C}_{\cal N}}(\langle t_{i},t_{j} \rangle ) \ \mbox{etc},
\end{displaymath}
single point and two point stabilizers in $\cal N$ respectively.
The coset stabilizing subgroup, ${\cal N}^{(w)}$, of $\cal N$ is
given by,
\begin{displaymath}
{{\cal N}^{(w)}} = \{{\pi} \in {\cal N}: {\cal N}w{\pi} = {\cal
N}w\},
\end{displaymath}
for $w$ a word in the symmetric generators. Clearly ${\cal N}^{w}
\leq {\cal N}^{(w)}$, and the number of cosets in the double coset
$[w] = {\cal N}w{\cal N}$ is given by $\mid {\cal N} \mid / \mid
{{\cal N}^{(w)}} \mid,$ since
\begin{eqnarray*}
{\cal N}w{\pi}_{1} \neq {\cal N}w{\pi}_{2} &\iff& {\cal
N}w{\pi}_{1}{{\pi}_{2}}^{-1} \neq {\cal N}w \cr &\iff&
{\pi}_{1}{{\pi}_{2}}^{-1} \notin {\cal N}^{(w)} \cr &\iff& {\cal
N}^{(w)}{\pi}_{1}{{\pi}_{2}}^{-1} \neq {\cal N}^{(w)} \cr &\iff&
{\cal N}^{(w)}{\pi}_{1} \neq {\cal N}^{(w)}{\pi}_{2}.
\end{eqnarray*}
In order to obtain the index of $\cal N$ in $\cal G$ we shall
perform a manual double coset enumeration of $\cal G$ over $\cal
N$; thus we must find all double cosets $[w]$ and work out how
many single cosets each of them contains. We shall know that we
have completed the double coset enumeration when the set of right
cosets obtained is closed under right multiplication. Moreover,
the completion test above is best performed by obtaining the
orbits of ${\cal N}^{(w)}$ on the symmetric generators. We need
only identify, for each $[w]$,  the double coset to which ${\cal
N}wt_{i}$ belongs for one symmetric generator $t_{i}$ from each
orbit. We will decompose the image ${\cal G}$ into double cosets
${\cal N} g {\cal N}$, where $g \in 2^{*n}:{\cal N}$ and find a
set $\{g_1, g_2, \ldots \}$ of elements of
${\cal G}$ such that\\
${\cal G} = {\cal N} g_1 {\cal N} \cup {\cal N} g_2 {\cal N}
\cup \ldots $. \\
 But for each $i$, we have $g_i = \pi_i w_i$, where $\pi_i \in
 {\cal N}$ and $w_i$ is a word in the  $t_{i}$, and so the double
 coset deomposition simplifies to \\
${\cal G} = {\cal N}\cup {\cal N} w_2 {\cal N} \cup {\cal N} w_3
{\cal N} \cup
\ldots $, \\
where $w_1$ is chosen to be the identity. When the set of
relations by which we are factoring is empty this gives the double
coset decomposition of the progenitor $2^{*n}:{\cal N}$, and in
this case there are infinitely many double cosets corresponding to
the orbits of ${\cal N}$ on the ordered $k$-tuples of the letters
of $\Omega = \{t_1, \ldots , t_n\}$ which have no adjacent
repetitions, where $k \in \mathbb{N}= \{1,2, \ldots \}$.
\par
We now give two interesting examples to illustrate the process.  \\
{\bf Example 1} \\
We consider
\begin{equation}
{\cal G} \cong
\frac{2^{*6}:L_{2}(5)}{[(\infty,0,1)(2,4,3)t_{2}]^{5}}
\newline \cong \langle {\cal N}, {\cal T} \mid {\cal N} \cong
L_{2}(5), t_{i}^{\pi} = t_{\pi (i)},
[(\infty,0,1)(2,4,3)t_{2}]]^{5} = 1 \rangle ,
\end{equation}
Thus $\cal G$ is the homomorphic image of the split extension
$2^{*6}:L_{2}(5)$ factored by the relator
$[(\infty,0,1)(2,4,3)t_{2}]^{5}$ and the action of ${\cal N} =
L_{2}(5)$ on the six symmetric generators is given by $x \sim
(0,1,2,3,4)$ and $y \sim (0,\infty)(1,4)$. Of course at this stage
the group  could still be infinite, or could collapse to the
identity. The given
relation\\
$[\pi t_{2}]^{5} = 1$, where $\pi = (\infty,0,1)(2,4,3)$, implies
\begin{eqnarray*}
\pi t_{2}\pi t_{2}\pi t_{2}\pi t_{2}\pi t_{2} = 1\\
\Rightarrow \pi^2 (\pi)^{-1}t_{2}\pi t_{2}\pi^3
(\pi)^{-2}t_{2}\pi^2 (\pi)^{-1}t_{2}\pi t_{2} = 1\\
\Rightarrow \pi^2 t_{2}^{\pi} t_{2}t_{2}^{\pi^2} t_{2}^{\pi} t_{2} = 1\\
\Rightarrow \pi^2 t_{4} t_{2}t_{3}t_{4} t_{2} = 1\\
\Rightarrow t_{4} t_{2}t_{3}t_{4} t_{2} = \pi.
\\
\end{eqnarray*}
Note that, in particular, all conjugates of the relator $\pi^2
t_{4} t_{2}t_{3}t_{4} t_{2}$, under conjugation by ${\cal N}$, are
also relators. We identify elements of the progenitor
$2^{*6}:L_{2}(5)$ with their images in ${\cal G}$ because in the
cases that interest us the image of ${\cal N}$ will be isomorphic
to ${\cal N}$ and the images of $t_{i}$ will be of order 2 and be
distinct. The vertices of the Cayley graph of $2^{\star n}: {\cal
N}$ has
the set of right cosets \\
$\{{\cal N}w$ $|$ w a word in the symmetric generators $t_i\}$\\
as its vertices. Each vertex ${\cal N}w$ is joined to the vertex
${\cal N}wt_i$, for $i \in \{1,2, \ldots , n\}$. The collapsed
Cayley graph is a diagram of the Cayley graph in which each orbit
of ${\cal N}$ in its action on the vertices by right
multiplication is represented by a single node, labelled with the
number of vertices that it contains. Lines between these nodes are
labelled with integers to indicate how many edges from a vertex of
one node lead to vertices of the other. The orbits of ${\cal N}$
on the right cosets of ${\cal N}$ in ${\cal G}$ are, of course,
just the double cosets of the form ${\cal N}g{\cal N}$ and the
labels on the nodes indicate how many single cosets a given double
coset contains. So in Figure the the node ${\cal
N}t_{\infty}t_0t_1$ is labelled 20, since ${\cal
N}t_{\infty}t_0t_1 = {\cal N}t_1t_2t_4 = {\cal
N}t_4t_3t_{\infty}$. Generally, Cayley graphs contain multiple
edges and loops. For example, the graph below does not have
multiple edges but has two loops, the integer 1+2 at the top of
the loop over the node labelled 30 indicates that $[\infty 0
\infty] = [\infty 0]$ and $[\infty 0 1] = [\infty 0]$, with
$\infty$ in the 1-orbit $\{\infty\}$ and 1 in the 2-orbit
$\{1,4\}$. The set of all double cosets [$w$]=${\cal N}w{\cal N}$,
coset stabilizing subgroups
${\cal N}^{(w)}$, and the number of single cosets each contains are exhibited in Table 1. \\ \\
\begin{figure}[h]
\begin{center}
\begin{picture}(265,60)
\put(18,28){\circle{30}} \put(17,25){1} \put(12,0){$[\star]$}
\put(32,28){\line(1,0){38}} \put(35,30){\small 6}
\put(63,30){\small 1} \put(84,28){\circle{30}} \put(79,25){6}
\put(78,0){$[\infty]$} \put(98,28){\line(1,0){38}}
\put(101,30){\small 5} \put(129,30){\small 1}
\put(150,28){\circle{30}} \put(140,25){30} \put(141,0){$[\infty 0
]$} \put(166,49){\circle{21}} \put(160,60){\small 1+2}
\put(164,28){\line(1,0){38}} \put(167,30){\small 2}
\put(194,30){\small 3} \put(216,28){\circle{30}} \put(207,25){20}
\put(204,0){$[\infty 0 1]$} \put(232,49){\circle{21}}
\put(231,60){\small 3}
\end{picture}
\end{center}
\vspace{.05in} \caption{The Cayley graph of $L_2(19)$ over
$L_2(5)$}
\end{figure}
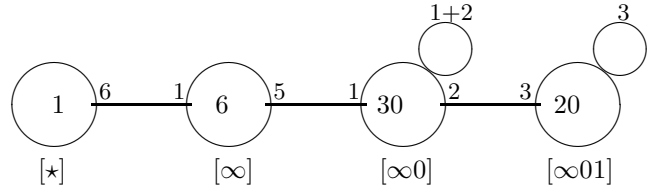
\par

\newpage
\begin{table}
\begin{tabular}{l l r}
Label [$w$] & Coset stabilizing subgroup ${\cal N}^{(w)}$ & Number
\cr
 & & of cosets  \cr \hline
\\
$[\star]$ & Since $\cal N$ is transitive on T & 1 \cr
\\
$[\infty]$ & \parbox[t]{4 in}{\raggedright ${\cal N}^{(\infty)} =
N^\infty =\langle x,(1,3)(0,4)\rangle \cong D_{10}$, has orbits
$\{\infty\}$ and $\{0,1,2,3,4\}$ on $T$.} & 6 \cr
\\
$[\infty0]$ & \parbox[t]{4 in}{\raggedright ${\cal N}^{\infty 0} =
N^{(\infty 0)}=\langle (1,4)(2,3) \rangle \cong C_2$, has orbits
$\{0\}$, $\{\infty\}$, $\{1, 4\}$ and
$\{2,3\}$.} & 30 \cr \\
$[\infty 0 0] = [\infty]$ & \parbox[t]{4 in}{\raggedright We show
that $\infty 0 \infty \sim \infty 0$.\\
Firstly, $\infty . 01401 = \infty . (1,4,0)(2,3,\infty )$\\
$\implies$$\infty . 01401 = (1,4,0)(2,3,\infty )2$\\
$\implies$$\infty . 014 = (1,4,0)(2,3,\infty )210$\\
Thus $\infty 014 \sim 210$ and $210 = (1,2,0)(3, \infty,4)12$
$\implies$ $\infty 014 \sim 12$. So $\infty 01 \sim 124$.\\
Secondly, $\infty . 04104 = \infty .(1,4,0)(2, \infty,3)$\\
$\implies$$\infty . 04104 = (1,4,0)(2, \infty,3)3$.\\
Thus $\infty 041 \sim 340 \sim 43$, since $34034 =
(1,2,\infty)(3,0,4)$ $\implies$$\infty 04 \sim 431$.\\
Now $\infty 0 \infty = \infty 0 1.1\infty \sim 124.1\infty$ =
$12.41\infty 41.14 = 12.(1,4,\infty)(2,3,0).14 \sim 4314 \sim
\infty 044 \sim \infty 0$.
}\\
$[\infty 0\infty] = [\infty 0]$ &\parbox[t]{4 in}{\raggedright and
$\infty 0 2\infty 0 = (1,3,4)(2,0,\infty)$ $\implies$ $\infty 0 2
\sim  0 \infty$} \\
$[\infty 02] =[\infty 0]$ &
 \parbox[t]{4 in}{\raggedright $\infty 0 1 = \infty 0 3.31
 \sim 0 \infty 31 \sim 0 (1,3,\infty)(2,0,4)3\infty)3\infty
 \sim 43\infty$. Thus $\infty 0 1 \sim 43\infty \sim 124$ } \\
$[\infty 01]$ & \parbox[t]{4 in}{\raggedright ${\cal N}^{(\infty
01)} =\langle (1,4,\infty )(2,3,0) \rangle \cong C_3$, has orbits
$\{1, 4,\infty \}$ and $\{2,3,0\}$ on ${\cal T}$.\\
Now $\infty 012 \sim \infty.(1,0,2)(3,4,\infty )10 \sim 310$.} & 20 \cr \\
$[\infty 012]=[\infty 01]$\\
$[\infty 011]= [\infty 0]$\\
\\ \\ \hline
\end{tabular}
\caption{The double cosets [$w$]=${\cal N}w{\cal N}$ in
$L_2(19)$.}
\end{table}

\par
\noindent
\\
Our argument shows that the maximum possible index of ${\cal N}$
in ${\cal G}$ is\\
$\frac{{\cal N}}{\cal N} + \frac{{\cal N}}{{\cal N}^{({\infty})}}
+ \frac{{\cal N}}{{\cal N}^{({\infty}0)}} + \frac{{\cal N}}{{\cal
N}^{({\infty}01)}} = 1 + \frac{{\cal N}}{|D_{10}|} + \frac{{\cal
N}}{|C_2|}+ \frac{{\cal N}}{|C_3|}= 1 + 6 + 30 + 20 = 57$. Thus
$|{\cal G}| \le 57 \times |N| = 57 \times 60 = 3420$. We now show
that G is isomorphic to $L_2(19)$. We construct a homomorphism
from the progenitor $2^{*6}:L_{2}(5)$ to $L_2(19)$ by defining\\
$x \equiv (\frac{8\eta + 15}{\eta})= (1, 4, 7, 2, 6)(3, 13, 15,
9, 16)(5, 11, 18, 12, 14)(8, 17, 10, 0, \infty)$, and \\
$y \equiv (\frac{-\eta - 1}{2\eta + 1})= (1, 12)(2, 7)(3, 13)(4,
10)(5, 15)(6, 17)(8, 14)(9, 20)(11, 16)(18, 0)$.\\
Since the order of $xy$ is 5, ${\cal N} = <x,y> \cong A_5$. We now
let\\ $t_{\infty} \equiv (\frac{4\eta + 15 }{\eta - 4})= (1, 0)(2,
17)(3, 11)(4, \infty)(5, 16)(6, 10)(7, 8)(9, 14)(12, 15)(13,
18)$,\\ and find that\\ $|t_{\infty}^{\cal N}| = 6$ and that
${\cal N}$
permutes the six images of $t_{\infty}$, by conjugation, as the group $L_2(5)$. Thus\\
$t_0= (1, 5)(2, 14)(3,19)(4, 17)(6, 7)(8, 20)(9, 10)(11, 15)(12,
18)(13,16)$,\\
$t_1= (1, 2)(3, 15)(4, 11)(5, 6)(7, 10)(8, 17)(9, 18)(12, 14)(13,
20)(16, 19)$,\\
 $t_2= (1, 11)(2, 19)(3, 20)(4, 6)(5, 14)(7, 18)(8,
15)(9, 13)(10, 17)(12, 16)$,\\
$ t_3= (1, 7)(2, 12)(3, 14)(4, 18)(5,
11)(6, 20)(8, 13)(9, 17)(10, 19)(15, 16)$,\\
 $t_4= (1, 8)(2, 4)(3, 9)(5,
13)(6, 14)(7, 12)(10, 16)(11, 18)(15, 17)(19, 20)$,\\
$t_{\infty}= (1, 19)(2, 17)(3, 11)(4, 20)(5, 16)(6, 10)(7, 8)(9,
14)(12, 15)(13, 18)$,\\
 where, under conjugation, we have:\\
 $x$: $(t_1,t_2,t_3,t_4,t_0)(t_{\infty})$ and \\
 $y$ : $(t_{\infty},t_0)(t_1,t_4)(t_2)(t_3)$.\\
 Since ${\cal N} \cong A_5$ is maximal in $L_2(19)$ and $t_{\infty}
 \notin {\cal N}$, we can conclude that $L_2(19)$ is a homomorphic
 image of the progenitor $2^{*6}:L_{2}(5)$; it remains to
 check that the additional relation $[\pi t_{2}]^{5} = 1$
 $\iff$ $t_{4} t_{2}t_{3}t_{4} t_{2} = \pi = (\infty,0,1)(2,4,3) = ((x*y)^{(x^2*y*x^3*y*x^3)})$
 $\iff$ $((x*y)^{(x^2*y*x^3*y*x^3)}*t_2)^5 = 1$ holds. We readily
 calculate that\\
 $t_{4} t_{2}t_{3}t_{4} t_{2} =(1, 17, 14)(3, 11, 15)(4, 12, 19)(5, 8, 9)(6, 7, 10)(16, 20,
 18)$ \\ = $((x*y)^{(x^2*y*x^3*y*x^3)}$,\\ which acts as
 $(\infty,0,1)(2,4,3)$, by conjugation, on the six symmetric
 generators. This shows that $L_2(19)$ is an image of $G$; but $|G|
 \le 57 \times 60 = |L_2(19)|$, and so the equality holds and $G
 \cong L_2(19)$.\\
The Cayley graph of $L_2(19)$ over $L_2(5)$ suggests that every
 element of $L_2(19)$ can be written as a permutation, of $L_2(5)$ on
 six letters followed by a word, in term of the symmetric
 generators, of length at most three.\\\\
We now give an example of a progenitor where the homomorphic image
$\cal G$ does not possess a symmetric generating set. In fact, the
symmetric generators will generate a subgroup of $\cal G$ of index
2.\\
{\bf Example 2} \\
We consider
\begin{equation}
{\cal G} \cong \frac{2^{*3}:S_3}{{[(0,1,2)t_0]^{10},[(0,1)t_0]^6}}
\newline \cong \langle {\cal N}, {\cal T} \mid {\cal N} \cong
S_3, t_{i}^{\pi} = t_{\pi (i)}, [(0,1,2)t_0]^{10} =1, [(0,1)t_0]^6
\rangle ,
\end{equation}
The group ${\cal G}$, given above, is isomorphic to $5^2 : D_6$.\\
{\bf Manual double coset enumeration of $5^2 : D_6$ over $S_3$}:\\
We first note that the double coset $[0]= {\cal N} t_0 {\cal N} =
\{{\cal N} t_0, {\cal N} t_1, {\cal N} t_2\}$ and the double coset
$[01] = {\cal N} t_0t_1 {\cal N} = \{{\cal N} t_0t_1, {\cal
N}t_0t_2, {\cal N}t_1t_0, {\cal N}t_1t_2, {\cal N}t_2t_0, {\cal
N}t_2t_0\}$. In order to compute the remaining double cosets, we
first examine our relations. The relation $[\pi t_0 ]^{10} = 1$,
where $\pi = (0,1,2)$,
gives\\
$\pi t_0\pi t_0\pi t_0\pi t_0 \pi t_0\pi t_0\pi t_0
\pi t_0\pi t_0\pi t_0 = 1$\\
 $\implies$ $\pi \pi \pi^{-1}t_0\pi t_0\pi \pi \pi^{-1}t_0\pi
 t_0\pi \pi \pi^{-1}t_0\pi t_0\pi \pi \pi^{-1}t_0
 \pi t_0\pi \pi \pi^{-1}t_0\pi t_0 = 1$\\
 $\implies$ $\pi^{-1}t_1t_0\pi^{-1}t_1t_0\pi^{-1}t_1t_0\pi^{-1}t_1t_0\pi^{-1}t_1t_0 =
 1$\\
  $\implies$
  $\pi^{-1}\pi^{-1}\pi t_1t_0\pi^{-1}t_1t_0\pi^{-1}\pi^{-1}\pi t_1t_0
  \pi^{-1}t_1t_0\pi^{-1}t_1t_0 = 1$\\ $\implies$\\
  $\pi t_0t_2t_1t_0t_2t_1t_0t_2t_1t_0 = 1$\\
  (Note that $[\pi t_0 ]^{10} = 1$ can be written as\\
$\pi^{10}t_0^{\pi^9}t_0^{\pi^8}t_0^{\pi^7}t_0^{\pi^6}t_0^{\pi^5}t_0^{\pi^4}
t_0^{\pi^3}t_0^{\pi^2}t_0^{\pi}t_0 = 1$).\\ Thus
$(0,1,2)t_0t_2t_1t_0t_2t_1t_0t_2t_1t_0 = 1$ and
\begin{eqnarray} (0,1,2)t_0t_2t_1t_0t_2 = t_0t_1t_2t_0t_1. \end{eqnarray}
Also, the relation $[(0,1)t_0]^6 = 1$ $\implies$
$(0,1)t_0(0,1)t_0(0,1)t_0(0,1)t_0(0,1)t_0(0,1)t_0 =
t_0^{(0,1)}t_0t_0^{(0,1)}t_0t_0^{(0,1)}t_0 = t_1t_0t_1t_0t_1t_0 =
1$ $\implies$ \begin{eqnarray} t_0t_1t_0 = t_1t_0t_1\end{eqnarray} .\\
Now $t_0t_1t_2t_0t_2t_1t_0 = t_0t_1.t_2t_0t_2.t_1t_0 =
t_0t_1.t_0t_2t_0.t_1t_0 = t_0t_1t_0.t_2.t_0t_1t_0 =
t_1t_0t_1.t_2.t_1t_0t_1 = t_1t_0.t_1t_2t_1.t_0t_1 =
t_1t_0.t_2t_1t_2.t_0t_1$.
This gives
\begin{eqnarray} t_0t_1t_2t_0t_2t_1t_0 = t_1t_0t_2t_1t_2t_0t_1
\end{eqnarray}.\\
Now ${\cal N}^{01} = N^{(01)}\ge < e >$, has orbits $\{0\}$,
$\{1\}$ and $\{ 2\}$. Thus we must consider the double cosets
$[010]$, $[011]$ and $[012]$. However,$[011] = [0]$ and $[012]=
{\cal N} t_0t_1t_2{\cal N} = \{{\cal N} t_0t_1t_2, {\cal
N}t_0t_2t_1, {\cal N}t_1t_0t_2, {\cal N}t_1t_2t_0, {\cal
N}t_2t_0t_1, {\cal N}t_2t_1t_0\}$. Moreover, (by )$010 = 101$
$\implies$ $|{\cal N}^{(010)}| \ge |<(0,1)>| = 2$. Thus the double
coset $[010]$ contains $\frac{|{\cal N}|}{|{\cal N}^{(010)}|} = 3$
single cosets. Now $N^{(010)}\ge <(0,1)> \cong C_2$ has orbits
$\{0,1\}$ and $\{2\}$ on ${\cal T}$ and $0102 \sim 1012$. Thus the
double coset $[0102]$ contains three single cosets. ${\cal
N}^{(0102)}= <(0,1)> \cong C_2$ has orbits $\{0,1\}$ and $\{2\}$
on ${\cal T}$ and $020 = 202 \implies 01020 = 01202$. Thus
$[01020]=[01202]$ and the double coset has six single cosets. At
this point we have the double cosets $[012]$ and $[01202]$ left to be considered.\\
Now $N^{(012)}= <e>$ has orbits $\{0\}$, $\{1\}$ and $\{2\}$ on
${\cal T}$. But $121 = 212 \implies 0121 = 0212 \implies$ the
coset $0121$ has two names. Thus $[0121]$ contains 3 single
cosets. Now $N^{(0121)}\ge <(1,2)> \cong C_2$ has orbits $\{0\}$
and $\{1,2\}$ on ${\cal T}$ and $01210 = 02120$ $\implies$ the
coset $01210$ has two names. Thus $[01210]$ contains three single
cosets and $N^{(01210)}\ge <(1,2)> \cong C_2$ has orbits $\{0\}$
and $\{1,2\}$ on ${\cal T}$ with $012101 = 012010$ $\implies$
$[012101]=[012010]$ and since $(0,1,2)t_0t_2t_1t_0t_2 =
t_0t_1t_2t_0t_1$ $\implies$$t_0t_1t_2t_0t_1t_0 =
(0,1,2)t_0t_2t_1t_0t_2t_0$, and therefore $012010 \sim 021020$,the
double coset $[012010]$ contains three single cosets. We now have
the double cosets  $[0120]$ and $[01202]$ that must be considered.
Now $N^{(0120)} \ge <e>$ has orbits $\{0\}$, $\{1\}$ and $\{2\}$
on ${\cal T}$. Because of the relation $01201 = (0,1,2)02102$, the
coset $01201$ has two names. Thus $[01201]$ contains three single
cosets. Since ${\cal N}^{(01201)} \ge <(1,2)> \cong C_2 $ has
orbits $\{0\}$, $\{1,2\}$ on ${\cal T}$ and the double coset
$[012010]$ has already been investigated. We are now left with the
double coset $[01202]$. From above, $[01202]$ cotains 6 single
cosets. Since $N^{(01202)} \ge <e> $ has orbits $\{0\}$, $\{1\}$
and $\{2\}$ on ${\cal T}$ and $01201 = (0,1,2)02102$ $\implies$
$0120 = (0,1,2)021021$ $\implies$ $012021 = (0,1,2)02102121$
$\implies$ $012021 = (0,1,2)021012$, the coset $012021$ has two
names. Thus the double coset $[012021]$ contains 3 single cosets.
$N^{(012021)} \ge <(1,2)> \cong C_2$ has orbits $\{0\}$ and
$\{1,2\}$ on ${\cal T}$. From above, $012021 = (0,1,2)021012$,
$0120210 =(0,1,2)0210120 = (0,1,2)0201020 = (0,1,2)2021020 =
(0,1,2)2021202 =(0,1,2)2012102$. Thus $0120210 =(0,1,2)0210120 =
(0,1,2)2012102$. Then, conjugation by $(0,1)$ gives, $1021201
=(1,0,2)1201021 = (1,0,2)2102012$ and $0120210 = 1021201$ (page ).
Hence $0120210 =(0,1,2)0210120 = (0,1,2)2012102 = 1021201
=(0,1,2)1201021 = (0,1,2)2102012$. $[0120210]$ contains 1 coset
and $N^{(0120210)} \ge <(1,2),(0,1,2)> \cong S_3$ is transitive on
${\cal T}$.
\newpage
\begin{table}
\begin{tabular}{l l r}
Label [$w$] & Coset stabilizing subgroup ${\cal N}^{(w)}$ & Number
\cr
 & & of cosets  \cr \hline
\\
$[\star]$ & Since $\cal N$ is transitive on T & 1 \cr
 $[0]$ &
\parbox[t]{4 in}{\raggedright ${\cal N}^{(0)} = N^0 \ge
\langle(1,2)\rangle \cong C_{2}$, has orbits $\{0\}$ and $\{1,2\}$
on ${\cal T}$.} & 6 \cr $[01]$ & \parbox[t]{4 in}{\raggedright
${\cal N}^{01} = N^{(01)}\ge \langle e \rangle $, has orbits
$\{0\}$, $\{1\}$ and $\{ 2\}$. Moreover, $010 = 101$ $\implies$
$|N^{(010)}| \ge |<(0,1)>| = 2$}. & 6 \cr $[010]$ &\parbox[t]{4
in}{\raggedright $N^{(010)}\ge <(0,1)> \cong C_2$ has orbits
$\{0,1\}$ and $\{2\}$ on ${\cal T}$ and $0102 \sim 1012$. Thus the
coset $0102$ has two names. This implies} & 3 \cr $[0102]$
&\parbox[t]{4 in}{\raggedright $N^{(0102)}= <(0,1)> \cong C_2$ has
orbits $\{0,1\}$ and $\{2\}$ on ${\cal T}$ and $020 = 202 \implies
01020 = 01202$.}  & 3 \cr $[01021]=[01202]$ & & 6\cr $[012]$ &
\parbox[t]{4 in}{\raggedright $N^{(012)}= <e>$ has orbits $\{0\}$,
$\{1\}$ and $\{2\}$ on ${\cal T}$.}& 6 \cr $[0120]$ &
\parbox[t]{4 in}{\raggedright Now $121 = 212 \implies 0121 =
0212 \implies$ the coset $0121$ has two names.}  & 6 \cr $[0121]$
&\parbox[t]{4 in}{\raggedright $N^{(0121)}\ge <(1,2)> \cong C_2$
has orbits $\{0\}$ and $\{1,2\}$ on ${\cal T}$ and $01210 = 02120$
$\implies$ the coset $01210$ has two names. }& 3 \cr $[01210]$
&\parbox[t]{4 in}{\raggedright $N^{(01210)}\ge <(1,2)> \cong C_2$
has orbits $\{0\}$ and $\{1,2\}$ on ${\cal T}$ and $012101 =
012010$ and the relation $01201 = (0,1,2)02102$ implies $012010 =
(0,1,2)021020$ and so the coset $012010$ has two names.} & 3 \cr
$[012101]=[012010]$ &
\parbox[t]{4 in}{\raggedright $N^{(0120)} \ge <e>$ has orbits
$\{0\}$, $\{1\}$ and $\{2\}$ on ${\cal T}$. Because of the
relation $01201 = (0,1,2)02102$, the coset $01201$ has two names.}
& 3\cr $[01201]$ &
\parbox[t]{4 in}{\raggedright $N^{(01201)} \ge <(1,2> \cong C_2 $
has orbits $\{0\}$, $\{1,2\}$ on ${\cal T}$ and the relation
$01201 = (0,1,2)02102$ implies $012010 = (0,1,2)021020$ and so the
coset $012010$ has two names.} & 3\cr $[012010]$ & \parbox[t]{4
in}{\raggedright $N^{(01202)} \ge <e> $ has orbits $\{0\}$,
$\{1\}$ and $\{2\}$ on ${\cal T}$ and the relation $01201 =
(0,1,2)02102$ $\implies$ $0120 = (0,1,2)021021$ $\implies$ $012021
= (0,1,2)02102121$ $\implies$ $012021 = (0,1,2)021012$ implies
that the coset $012021$ has two names.} & 3\cr $[012021]$ &
\parbox[t]{4 in}{\raggedright $N^{(012021)} \ge <(1,2)> \cong C_2$
has orbits $\{0\}$ and $\{1,2\}$ on ${\cal T}$. From above,
$012021 = (0,1,2)021012$, $0120210 =(0,1,2)0210120 =
(0,1,2)0201020 = (0,1,2)2021020 = (0,1,2)2021202 =(0,1,2)2012102$.
Thus $0120210 =(0,1,2)0210120 = (0,1,2)2012102$. Then, conjugation
by $(0,1)$ gives, $1021201 =(1,0,2)1201021 = (1,0,2)2102012$ and
$0120210 = 1021201$ (page ). Hence $0120210 =(0,1,2)0210120 =
(0,1,2)2012102 = 1021201 =(0,1,2)1201021 = (0,1,2)2102012$. } &
3\cr $[0120210]$ &
\parbox[t]{4 in}{\raggedright $N^{(0120210)} \ge <(1,2),(0,1,2)>
\cong S_3$ is transitive on  ${\cal T}$. } & 1\cr \hline
\end{tabular}
\caption{The double cosets [$w$]=${\cal N}w{\cal N}$ in
$5^2:D_6$.}
\newpage
\end{table}
\begin{figure}
\setlength{\unitlength}{.65pt}
\begin{picture}(300,150)(1,1)
\put(18,28){\circle{30}} \put(17,25){1} \put(12,0){$[\star]$}
\put(32,28){\line(1,0){38}} \put(36,30){\small 3}
\put(63,30){\small 1} \put(84,28){\circle{30}} \put(79,25){3}
\put(78,0){$[0]$} \put(98,28){\line(1,0){38}} \put(102,30){\small
2}\put(127,30){\small 1}\put(166,30){\small 1}
\put(150,28){\circle{30}}\put(145,25){6} \put(141,0){$[01]$}
\put(201,79){\line(-1,-1){39}}\put(158,44){\small 1}
\put(193,79){\small 2} \put(215,79){\circle{30}} \put(213,74){3}
\put(205,52){$[010]$} \put(230,79){\line(1,0){38}}
\put(233,81){\small 1}\put(282,79){\circle{30}} \put(279,76){3}
\put(268,52){$[0102]$}\put(300,79){\small 2} \put(260,81) {\small
1} \put(297,79){\line(1,-1){40}} \put(333,45){\small 1}
\put(165,29){\line(1,0){38}}\put(217,29){\circle{30}}\put(215,25){6}
\put(194,30){\small 1}\put(233,30){\small
1}\put(179,2){$[012]$}\put(220,3){\small 1}
 \put(217,15){\line(0,-1){38}}
\put(217,-37){\circle{30}}\put(215,-39){3}\put(220,-19){\small 2}
\put(174,-64){$[0121]$} \put(220,-62){\small 1}
\put(217,-51){\line(0,-1){38}} \put(217,-103){\circle{30}}
\put(220,-85){\small 1} \put(235,-102){\small 2}
 \put(215,-109){3}
\put(195,-132){$[01210]$} \put(231,-103){\line(1,0){38}}
\put(260,-102){\small 2} \put(283,-103){\circle{30}}
 \put(281,-105){3}\put(260,-129){$[012010]$}
\put(276,-85){\small 1} \put(276,-63){\small
1}\put(290,-64){$[01201]$} \put(283,-89){\line(0,1){38}}
\put(283,-37){\circle{30}} \put(281,-39){3}\put(276,-20){\small 2}
\put(283,-23){\line(0,1){38}}\put(276,3){\small 1}
 \put(290,2){$[0120]$}
\put(232,29){\line(1,0){38}}\put(284,29){\circle{30}}\put(282,27){6}
\put(262,30){\small 1}\put(299,31){\small 1}
\put(298,29){\line(1,0){38}}\put(328,31){\small 1}
\put(350,29){\circle{30}}\put(348,27){6} \put(335,0){$[01202]$}
\put(364,29){\line(1,0){38}}\put(367,30){\small 1}
\put(393,30){\small 2} \put(416,29){\circle{30}}\put(414,28){3}
\put(399,0){$[012021]$} \put(430,29){\line(1,0){38}}
\put(433,30){\small 1} \put(459,30){\small
3}\put(482,29){\circle{30}}\put(480,28){1}
\put(462,0){$[0120210]$}
\end{picture}
\vspace{1.5in} \caption{The Cayley graph of $5^2:D_6$ over $S_3$}
\end{figure}
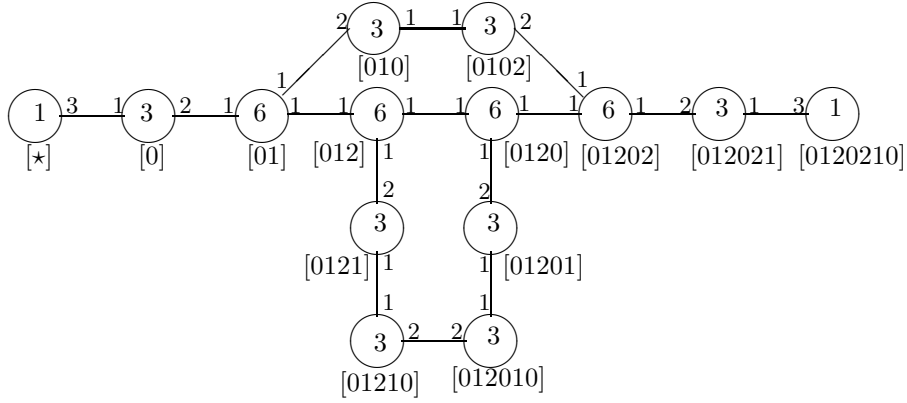
We now present our main problem. A presentation for the progenitor
$2^{*(7+7)}:PGL_2(7)$, is given by :\\
$<x,y,t|x^7, y^2, t^2, (x^{-1} * t)^2 , (y * x)^3 ,  t * x^{-1}*
y$\\
$* x * t * y , x^2 * y * x^3 * y * x^{-4} * y * x^{-4} * y * x,
s^2, (s^3,y),(s^4, x*y)>$,\\
where $N = PGL_2(7)  \cong
<x,y,t|x^7, y^2, t^2, (x^{-1} * t)^2 , (y
* x)^3 ,$\\
$t * x^{-1} * y* x * t * y , x^2 * y * x^3 * y * x^{-4} * y *
x^{-4} * y * x>$;\\
and the action of N on the symmetric generators is given by\\
       $x \sim  ({\bf 0},{\bf 1},{\bf 2},{\bf 3},{\bf 4},{\bf 5},{\bf 6})(0,6,5,4,3,2,1)$,\\
       $y \sim  ({\bf 2},{\bf 6})({\bf 4},{\bf 5})(0,3)(5,6)$,\\
       $t \sim ({\bf 0}, 0)({\bf 1},1)({\bf 2},2)({\bf 3},3)({\bf 4},4)({\bf 5},5)({\bf 6},6)$.\\
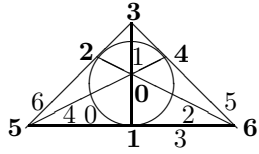
\begin{figure}
\setlength{\unitlength}{1pt}
\begin{center}
\begin{picture}(40,50)(75,166)
\put(97,200) {{\bf{3}}} \put(99,199){\line(-1,-1){40}}
\put(99,199){\line(1,-1){40}}\put(61,166){6}
\put(52,156){{\bf{5}}}\put(60,160){\line(1,0){78}}
\put(134,166){5}\put(141,156){{\bf{6}}} \put(81,161){0}
\put(115,152){3} \put(73,161){4}\put(118,161){2}
\put(97,152){{\bf{1}}} \put(99,176){\circle{30}}
\put(99,199){\line(0,-1){39}} \put(138,160){\line(-2,1){52}}
\put(79,185) {{\bf{2}}}\put(60,160){\line(2,1){52}}
\put(115,184){{\bf{4}}} \put(100,169){{\bf{0}}} \put(99,183){1}
\end{picture}
\end{center}
\caption{7-point Projective Plane.}
\end{figure}

We factor the progenitor by the following relations\\
$t= {\bf  0} 0 {\bf 0}$ and $y={\bf 1}  1 {\bf 1} 1$\\
to get the finite homomorphic image  G. Thus
\begin{equation}
G \cong \frac{2^{*(7+7)}:PGL_2(7)}{t= {\bf  0} 0 {\bf 0}, y=({\bf
0} 1)^2}.
\end{equation}
 \\
The index of $PGL_2(7)$ in G is 36. It turns out that $G \cong
U3(3):2$ (see \cite{JABBAR}).

\section{Manual double coset enumeration of $U_3(3):2$ over $PGL_2(7)$}
$t= {\bf 0} 0 {\bf 0} \implies {\bf 0} 0 = t {\bf 0} \implies {\bf
0} 0 \sim  {\bf 0}$.\\
$y=({\bf 0} 1)^2 = {\bf 0} 1{\bf 0} 1\implies {\bf 0} 1 \sim 1{\bf
0}$. Define $\pi_{{\bf 0}0{\bf 0}} = t = ({\bf 0}, 0)({\bf 1},1)({\bf 2},2)({\bf 3},3)({\bf 4},4)({\bf 5},5)({\bf 6},6)$,\\
and $\sigma_{{\bf 0}1} = {\bf 0} 1{\bf 0} 1 = y = ({\bf 2},{\bf 6})({\bf 4},{\bf 5})(0,3)(5,6)$.\\
We note that\\
$\pi_{{\bf 0}5{\bf 0}} = ({\bf 1}, 6)({\bf 2},3)({\bf 3},1)({\bf
4},4)({\bf 5},0)({\bf 6},2)({\bf 0},5)$ and \\
$\pi_{5{\bf 1}5} = ({\bf 1}, 5)({\bf 2},4)({\bf 3},1)({\bf
4},0)({\bf 5},2)({\bf 6},3)({\bf 0},6)$. Then\\
${\bf 0}{\bf 1} = {\bf 0}55{\bf 1}$\\
$ = {\bf 0}5.5{\bf 1}5.5 $\\
$={\bf 0}5{\bf 0}.{\bf 0}.5{\bf 1}5.5 $\\
$=({\bf 1}, 6)({\bf 2},3)({\bf 3},1)({\bf 4},4)({\bf 5},0)({\bf
6},2)({\bf 0},5).{\bf 0}.({\bf 1}, 5)({\bf 2},4)({\bf 3},1)({\bf
4},0)({\bf 5},2)({\bf 6},3)({\bf 0},6).5 $\\
$= {\bf 0}5.5{\bf 1}5.5 $\\
$={\bf 0}5{\bf 0}.{\bf 0}.5{\bf 1}5.5$\\
$ =({\bf 1}, 6)({\bf 2},3)({\bf 3},1)({\bf 4},4)({\bf 5},0)({\bf
6},2)({\bf 0},5)({\bf 1}, 5)({\bf 2},4)({\bf 3},1)({\bf
4},0)({\bf 5},2)({\bf 6},3)({\bf 0},6).65$.\\
Thus ${\bf 0}{\bf 1} = ({\bf 0},{\bf 1})({\bf 2},{\bf 6}), ({\bf
5},{\bf 4})(0,2,3,4)(5,6).65$.\\
So $[{\bf 0}{\bf 1}]^2 ={\bf 0}{\bf 1}{\bf 0}{\bf 1}$\\
$=({\bf 0},{\bf 1})({\bf 2},{\bf 6}, {\bf 5},{\bf
4})(0,2,3,4)(5,6).65.({\bf 0},{\bf 1})({\bf 2},{\bf 6}, {\bf
5},{\bf 4})(0,2,3,4)(5,6).65$\\
$= ({\bf 2},{\bf 5})({\bf 4},{\bf 6})(0,3)(2,4)$.\\
$\implies$ ${\bf 1}{\bf 0}=({\bf 2},{\bf 5})({\bf 4},{\bf
6})(0,3)(2,4){\bf 0}{\bf 1}$.\\
Define $\tau_{{\bf 0}{\bf 1}}= {\bf 0}{\bf 1}{\bf 0}{\bf 1}=
({\bf 2},{\bf 5})({\bf 4},{\bf 6})(0,3)(2,4)$.\\
Similarly,\\
${\bf 1}{\bf 0} = {\bf 1}66{\bf 0}$\\
$ = {\bf 1}6.6{\bf 0}6.6 $\\
$={\bf 1}6{\bf 1}.{\bf 1}.6{\bf 0}6.6 $\\
$=({\bf 1}, 6)({\bf 2},0)({\bf 3},1)({\bf 4},2)({\bf 5},3)({\bf
6},4)({\bf 0},5).{\bf 1}.({\bf 1}, 5)({\bf 2},2)({\bf 3},1)({\bf
4},3)({\bf 5},4)({\bf 6},0)({\bf 0},6).6 $\\
$= ({\bf 1}, 6)({\bf 2},0)({\bf 3},1)({\bf 4},2)({\bf 5},3)({\bf
6},4)({\bf 0},5)({\bf 1}, 5)({\bf 2},2)({\bf 3},1)({\bf
4},3)({\bf 5},4)({\bf 6},0)({\bf 0},6).56$\\
$=({\bf 0}, 1)({\bf 2},{\bf 6},{\bf 5},{\bf
4})(0,2,3,4)(5,6).56$.\\
Now $1 \sim {\bf 1}{\bf 0}0{\bf 0} \sim {\bf 0}{\bf 1}0{\bf 0}
\sim {\bf 0}.{\bf 1}0{\bf 1}0.0{\bf 1}{\bf 0} \sim {\bf 0}({\bf
3},{\bf 0})({\bf 5},{\bf 6})(2,6)(4,5)0{\bf 1}{\bf 0} \sim {\bf
3}0{\bf 1}{\bf 0} \sim {\bf 3}0({\bf 2},{\bf 5})({\bf 4},{\bf
6})(0,3)(2,4){\bf 0}{\bf 1} \sim {\bf 3}3{\bf 0}{\bf 1} \sim {\bf
3}{\bf 0}{\bf 1}$.\\
$\implies$ $1{\bf 1} \sim {\bf 3}{\bf 0}$. Hence $1{\bf 1} \sim
{\bf 3}{\bf 0} \sim {\bf 0}{\bf 3} \sim 03 \sim 30$.\\
$\pi_{{\bf 0}0{\bf 0}}1={\bf 1}\pi_{{\bf 0}0{\bf 0}}$\\
$={\bf 1}{\bf 0}0{\bf 0}= \tau_{{\bf 0}{\bf 1}}{\bf 0}{\bf 1}0{\bf
0}$\\
$= \tau_{{\bf 0}{\bf 1}}{\bf 0}\sigma_{{\bf 1}0}0{\bf 1}{\bf 0}$\\
$=\tau_{{\bf 0}{\bf 1}}\sigma_{{\bf 1}0}{\bf 3}0{\bf 1}{\bf 0}$\\
$=\tau_{{\bf 0}{\bf 1}}\sigma_{{\bf 1}0}{\bf 3}0\tau_{{\bf 0}{\bf
1}}{\bf 0}{\bf 1}$\\
$=\tau_{{\bf 0}{\bf 1}}\sigma_{{\bf 1}0}\pi_{{\bf 3}0{\bf 3}}{\bf
3}{\bf 1}{\bf 0} =\tau_{{\bf 0}{\bf 1}}\sigma_{{\bf 1}0}{\bf
3}0\tau_{{\bf 0}{\bf 1}}{\bf 0}{\bf 1}$\\
$=\tau_{{\bf 0}{\bf 1}}\sigma_{{\bf 1}0}\pi_{{\bf 3}0{\bf 3}}{\bf
3}\tau_{{\bf 0}{\bf 1}}{\bf 0}{\bf 1}$\\
$=\tau_{{\bf 0}{\bf 1}}\sigma_{{\bf 1}0}{\bf 3}0\tau_{{\bf 0}{\bf
1}}{\bf 0}{\bf 1} =\tau_{{\bf 0}{\bf 1}}\sigma_{{\bf 1}0}\pi_{{\bf
3}0{\bf 3}}\tau_{{\bf 0}{\bf 1}}{\bf 3}{\bf 0}{\bf 1}$\\
$=({\bf 2},{\bf 5})({\bf 4},{\bf 6})(0,3)(2,4).({\bf 0},{\bf
3})({\bf 5},{\bf 6})(6,2)(4,5).$\\
$({\bf 0},3)({\bf 1},1)({\bf 2},6)({\bf 3},0)({\bf 4},5)({\bf
5},2)({\bf 6},4).({\bf 2},{\bf 5})({\bf 4},{\bf
6})(0,3)(2,4){\bf 3}{\bf 0}{\bf 1}$.\\
$= ({\bf 0},3,{\bf 3},0)({\bf 1},1)({\bf 2},2,{\bf 6}, 5,{\bf
4},4,{\bf 5},6){\bf 3}{\bf 0}{\bf 1}$.\\
Thus $1 = ({\bf 2},{\bf 6})({\bf 4},{\bf 5})(0,3)(5,6){\bf 3}{\bf
0}{\bf 1}$.\\\\\\
\section{Manual double coset enumeration of $U_3(3):2$ over $PGL_2(7)$}
We first note that ${\cal N}t_{\bf{0}} {\cal N}$ contains fourteen
single cosets, namely, ${\cal N}t_{\bf{0}}, {\cal N}t_{\bf{1}},
\ldots , {\cal N}t_{\bf{6}}, {\cal N}t_0,{\cal N}t_1, \ldots,
{\cal N}t_6$. In order to compute the remaining double coset
${\cal N} w {\cal N}$, where $w$ is a word in the $t_is$, we first
analyze our relations.\\
The given relation $t= {\bf 0} 0 {\bf 0}$ gives ${\bf 0} 0 = t
{\bf 0}$. Thus ${\bf 0} 0 \sim  {\bf 0}$.\\
Now $y=({\bf 0} 1)^2 = {\bf 0} 1{\bf 0} 1\implies
{\bf 0} 1 \sim 1{\bf 0}$.\\
Define $\pi_{{\bf 0}0} = t = ({\bf 0}, 0)({\bf 1},1)({\bf 2},2)({\bf 3},3)({\bf 4},4)({\bf 5},5)({\bf 6},6)$,\\
and ${\bf 0} 1{\bf 0} 1 = y = ({\bf 2},{\bf 6})({\bf 4},{\bf
5})(0,3)(5,6)$ $\implies$ $[{\bf 0} 1{\bf 0} 1]^a = y^a$. Thus
\begin{eqnarray} {\bf 1} 0 {\bf 1} 0 = ({\bf 3, 0})({\bf 5, 6})(2,6)(4,5). \end{eqnarray}
 Now\\
${\bf 0}{\bf 1} = {\bf 0}55{\bf 1}$\\
$ = {\bf 0}5.5{\bf 1}5.5 $\\
$={\bf 0}5{\bf 0}.{\bf 0}.5{\bf 1}5.5 $\\
$= \pi_{{\bf 0}5}. {\bf 0}.\pi_{{\bf 1}5}.5$\\
$= \pi_{{\bf 0}5}.\pi_{{\bf 1}5}.{\bf 0}^{\pi_{{\bf 1}5}}.5$\\
 $=({\bf 1}, 6)({\bf 2},3)({\bf 3},1)({\bf 4},4)({\bf 5},0)({\bf
6},2)({\bf 0},5).({\bf 1}, 5)({\bf 2},4)({\bf 3},1)({\bf
4},0)({\bf 5},2)({\bf 6},3)({\bf 0},6).65 $\\
Thus ${\bf 0}{\bf 1} = ({\bf 0},{\bf 1})({\bf 2},{\bf 6,
5, 4})(0,2,3,4)(5,6).65$.\\
Similarly,\\
${\bf 1}{\bf 0} = {\bf 1}66{\bf 0}$\\
$ = {\bf 1}6.6{\bf 0}6.6 $\\
$={\bf 1}6{\bf 1}.{\bf 1}.6{\bf 0}6.6 $\\
$= \pi_{{\bf 1}6}.{\bf 1}.\pi_{{\bf 0}6}.6$\\
$= \pi_{{\bf 1}6}.\pi_{{\bf 0}6}.{\bf 1}^{\pi_{{\bf 0}6}}.6$\\
$= ({\bf 1}, 6)({\bf 2},0)({\bf 3},1)({\bf 4},2)({\bf 5},3)({\bf
6},4)({\bf 0},5)({\bf 1}, 5)({\bf 2},2)({\bf 3},1)({\bf
4},3)({\bf 5},4)({\bf 6},0)({\bf 0},6).56$\\
$=({\bf 0, 1})({\bf 2},{\bf 6},{\bf 5},{\bf
4})(0,2,3,4)(5,6).56$.\\
Thus ${\cal N}t_{{\bf 0}}t_{{\bf 1}} = {\cal N}t_6t_5 = {\cal
N}t_65t_6$.\\
Moreover, $[{\bf 0}{\bf 1}]^2 ={\bf 0}{\bf 1}{\bf 0}{\bf 1}$\\
$=({\bf 0},{\bf 1})({\bf 2},{\bf 6}, {\bf 5},{\bf
4})(0,2,3,4)(5,6).65.({\bf 0},{\bf 1})({\bf 2},{\bf 6}, {\bf
5},{\bf 4})(0,2,3,4)(5,6).65$\\
$= ({\bf 0},{\bf 1})({\bf 2},{\bf 6}, {\bf 5},{\bf
4})(0,2,3,4)(5,6).({\bf 0},{\bf 1})({\bf 2},{\bf 6}, {\bf 5},{\bf
4})(0,2,3,4)(5,6). [65]^{({\bf 0},{\bf 1})({\bf 2},{\bf 6}, {\bf
5},{\bf 4})(0,2,3,4)(5,6)}.65$\\
$= ({\bf 2},{\bf 5})({\bf 4},{\bf 6})(0,3)(2,4)56.65$\\
$= ({\bf 2},{\bf 5})({\bf 4},{\bf 6})(0,3)(2,4)$. Thus
\begin{eqnarray} [{\bf 0}{\bf 1}]^2 = ({\bf 2},{\bf 5})({\bf 4},{\bf 6})(0,3)(2,4). \end{eqnarray}
Also, \begin{eqnarray} {\bf 1}{\bf 0}=({\bf 2},{\bf 5})({\bf
4},{\bf
6})(0,3)(2,4){\bf 0}{\bf 1}.\end{eqnarray}\\
Hence \begin{eqnarray}{\bf 0}{\bf 1} \sim {\bf 1}{\bf 0} \sim 56
\sim 65.
\end{eqnarray}
 From (8) and (9) above, define $\sigma_{{\bf 0}1} = {\bf 0}
1{\bf 0} 1 = y = ({\bf 2},{\bf 6})
({\bf 4},{\bf 5})(0,3)(5,6)$ and\\
$\tau_{{\bf 0}{\bf 1}}= {\bf 0}{\bf 1}{\bf 0}{\bf 1}=
({\bf 2},{\bf 5})({\bf 4},{\bf 6})(0,3)(2,4)$.\\
Now \\
$\pi_{{\bf 0}0}1={\bf 1}\pi_{{\bf 0}0}$\\
$={\bf 1}{\bf 0}0{\bf 0}= {\bf 1}{\bf 0}.0{\bf 0} = \tau_{{\bf
0}{\bf 1}}{\bf 0}{\bf 1}.0{\bf
0}$\\
$= \tau_{{\bf 0}{\bf 1}}{\bf 0}.{\bf 1}0.{\bf
0}$\\
$= \tau_{{\bf 0}{\bf 1}}{\bf 0}\sigma_{{\bf 1}0}0{\bf 1}{\bf 0}$\\
$= \tau_{{\bf 0}{\bf 1}}{\bf 0}\sigma_{{\bf 1}0}0^{\sigma_{{\bf 1}0}}{\bf 1}{\bf 0}$\\
$= \tau_{{\bf 0}{\bf 1}}\sigma_{{\bf 1}0}{\bf 0}^{\sigma_{{\bf 1}0}}0{\bf 1}{\bf 0}$\\
$=\tau_{{\bf 0}{\bf 1}}\sigma_{{\bf 1}0}{\bf 3}0{\bf 1}{\bf 0}$\\
$=\tau_{{\bf 0}{\bf 1}}\sigma_{{\bf 1}0}\pi_{{\bf 3}0}{\bf 3}{\bf 1}{\bf 0}$\\
$=\tau_{{\bf 0}{\bf 1}}\sigma_{{\bf 1}0}\pi_{{\bf 3}0}{\bf
3}\tau_{{\bf 0}{\bf 1}}{\bf 0}{\bf 1}$\\
$=\tau_{{\bf 0}{\bf 1}}\sigma_{{\bf 1}0}\pi_{{\bf 3}0}\tau_{{\bf
0}{\bf 1}}{\bf
3}^{\tau_{{\bf 0}{\bf 1}}}{\bf 0}{\bf 1}$\\
$=\tau_{{\bf 0}{\bf 1}}\sigma_{{\bf 1}0}\pi_{{\bf 3}0}\tau_{{\bf
0}{\bf 1}} {\bf 3}{\bf
0}{\bf 1}$\\
Hence $\pi_{{\bf 0}0}1 = \tau_{{\bf 0}{\bf 1}}\sigma_{{\bf
1}0}\pi_{{\bf 3}0}\tau_{{\bf 0}{\bf 1}} {\bf 3}{\bf 0}{\bf 1}$.\\
So $1 = \pi^{-1}_{{\bf 0}0}\tau_{{\bf 0}{\bf 1}}\sigma_{{\bf
1}0}\pi_{{\bf 3}0}\tau_{{\bf 0}{\bf 1}} {\bf 3}{\bf 0}{\bf 1}$\\
$=({\bf 1}, 1)({\bf 2}, 2)({\bf 3}, 3)({\bf 4}, 4)({\bf 5},
5)({\bf 6}, 6)({\bf 7}, 7).({\bf 2},{\bf 5})({\bf 4},{\bf
6})(0,3)(2,4).({\bf 0},{\bf
3})({\bf 5},{\bf 6})(6,2)(4,5).$\\
$({\bf 0},3)({\bf 1},1)({\bf 2},6)({\bf 3},0)({\bf 4},5)({\bf
5},2)({\bf 6},4).({\bf 2},{\bf 5})({\bf 4},{\bf
6})(0,3)(2,4){\bf 3}{\bf 0}{\bf 1}$.\\
$= ({\bf 1}, 1)({\bf 2}, 2)({\bf 3}, 3)({\bf 4}, 4)({\bf 5},
5)({\bf 6}, 6)({\bf 7}, 7).({\bf 0},3,{\bf 3},0)({\bf 1},1)({\bf
2},2,{\bf 6}, 5,{\bf
4},4,{\bf 5},6){\bf 3}{\bf 0}{\bf 1}$.\\
Thus $1 = ({\bf 2},{\bf 6})({\bf 4},{\bf 5})(0,3)(5,6){\bf 3}{\bf
0}{\bf 1}$. Hence
\begin{eqnarray}{\bf 0}{\bf 3} =  ({\bf 2},{\bf 6})({\bf 4},{\bf 5})(0,3)(5,6){\bf 1}1. \end{eqnarray}
\\
We note that $1^{{\bf (1, 3, 0)(4, 5, 6)}(2, 0, 6)(3, 5, 4)} =$\\
$(({\bf 2},{\bf 6})({\bf 4},{\bf 5})(0,3)(5,6){\bf 3}{\bf 0}{\bf
1})^{{\bf (1, 3, 0)(4, 5, 6)}(2, 0, 6)(3, 5, 4)}$.\\
So \begin{equation} 1 =  {\bf (1, 3, 0)(2, 4, 6)}(2, 0, 5)(3, 6,
4){\bf 013}. \end{equation}
Now ${\cal N}^{({\bf 0})} = {\cal
N}^{\bf 0} = <({\bf 1, 2, 3, 6)(4, 5})(1, 2)(3, 0, 6, 5),{\bf (2,
4)(5, 6)}(2, 4)(5, 6)> \cong S_4$ has orbits $\{ {\bf 0} \}, \{1,
2, 4\}, \{3, 5, 6, 0 \}$, and $\{{\bf 1, 2, 3, 4, 5, 6}\}$. Now
the 1-orbit takes the double coset $[{\bf 0}]$ back to $[\star]$,
the 4-orbit loops back to itself, since ${\bf 0}0 = t {\bf 0}$ and
we need to examine the double cosets $[{\bf 0}1]$ and $[{\bf
01}]$. Now ${\bf 10} \sim {\bf 01}$ (from (9)) and form (10),
$[{\bf 0}{\bf 3}]^{{\bf{(1, 3)(4, 5)}}(3, 5)(6, 0)} = [({\bf
2},{\bf 6})({\bf 4},{\bf 5})(0,3)(5,6){\bf 1}1]^{{\bf{(1, 3)(4,
5)}}(3, 5)(6, 0)}$ gives
\begin{eqnarray} {\bf 0 1} = ({\bf 2},{\bf 6})({\bf 4},{\bf
5})(0,3)(5,6){\bf 3}1.
\end{eqnarray} Now from (9), $ (({\bf 2},{\bf 6})({\bf 4},{\bf
5})(0,3)(5,6){\bf 3}1)^2=({\bf 2},{\bf 5})({\bf 4},{\bf
6})(0,3)(2,4)$. \\
So $({\bf 2},{\bf 6})({\bf 4},{\bf 5})(0,3)(5,6){\bf 3}1.({\bf
2},{\bf 6})({\bf 4},{\bf 5})(0,3)(5,6){\bf 3}1 = ({\bf 2},{\bf
5})({\bf 4},{\bf 6})(0,3)(2,4)$ $\implies$ \begin{eqnarray} {\bf
3}1{\bf 3}1 =({\bf 2},{\bf 5})({\bf 4},{\bf 6})(0,3)(2,4).
\end{eqnarray} Thus ${\bf 3}1=({\bf 2},{\bf 5})({\bf 4},{\bf
6})(0,3)(2,4)1{\bf 3}$. Hence from (12), (14) and (15) , ${\bf
10}\sim {\bf 01} \sim {\bf 3}1 \sim 1{\bf 3} \sim 56 \sim 65$. Now
Now ${\cal N}^{{\bf 01}} = <y,(x * y * x^2)^2> = $\\ $<{\bf (2,
6)(4, 5)}(3, 0)(5, 6), {\bf (2, 5)(4, 6)}(2, 4)(3, 0)>$\\ and
    ${\cal N}t_{{\bf 0}}t_{{\bf 1}}{\bf (1}, 6){\bf (2}, 0){\bf (3}, 1)
    {\bf (4}, 2){\bf (5,} 3){\bf (6,} 4){\bf (0,} 5)$
    = ${\cal N}t_5t_6$ = ${\cal N}t_{\bf 0}t_{\bf 1}$.\\
Hence ${\cal N}^{({\bf 01})} \ge <{\bf (2, 6)(4, 5)}(3, 0)(5, 6),
    {\bf (2, 5)(4, 6)}(2, 4)(3, 0), {\bf (1}, 6){\bf (2}, 0){\bf (3}, 1)
    {\bf (4}, 2){\bf (5}, 3){\bf (6,} 4){\bf (0,} 5)> \cong
    D_{16}$.\\
Now the orbits of ${\cal N}^{({\bf 01})}$ on $\Omega$ are $\{{\bf
3}, 1 \}$, $\{ {\bf 0, 1}, 5, 6 \}$, and $\{ {\bf 2, 4, 5, 6}, 2,
3, 4, 0 \}$. Since ${\bf 013} \sim 1$, ${\bf 011 \sim 0}$, and
${\bf 01}2 = {\bf 0}.{\bf 1}2 = {\bf 0}\pi_{{\bf 1}2}{\bf 1} =$\\
$\pi_{{\bf 1}2}{\bf 0}^{\pi_{{\bf 1}2}}{\bf 1} = {\bf (1}, 2){\bf
(2}, 0){\bf (3}, 4){\bf (4}, 6){\bf (5}, 5){\bf (6}, 3){\bf (0},
1){\bf 0}^{{\bf (1}, 2){\bf (2}, 0){\bf (3}, 4){\bf (4}, 6){\bf
(5}, 5){\bf (6}, 3){\bf (0}, 1)}{\bf 1} =$\\
$ {\bf (1}, 2){\bf (2}, 0){\bf (3}, 4){\bf (4}, 6){\bf (5}, 5){\bf
(6}, 3){\bf (0}, 1)1{\bf 1} \sim {\bf 03} \in [{\bf 01}]$, we
must, therefore, have found all double cosets and the coset
diagram
indicated has form:\\\\\\
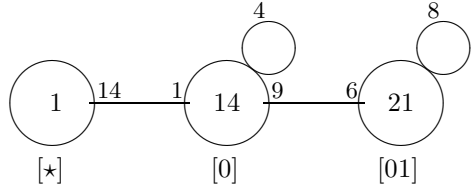
\begin{figure}
\begin{center}
\begin{picture}(200,50)
\put(18,28){\circle{30}} \put(17,25){1} \put(12,0){$[\star]$}
\put(32,28){\line(1,0){38}} \put(35,30){\small 14}
\put(63,30){\small 1} \put(84,28){\circle{30}} \put(79,25){14}
\put(78,0){$[0]$} \put(100,49){\circle{21}}\put(94,60){\small
4}\put(98,28){\line(1,0){38}} \put(101,30){\small
9}\put(129,30){\small 6}\put(150,28){\circle{30}}\put(145,25){21}
\put(141,0){$[01]$}\put(166,49){\circle{21}}\put(160,60){\small 8}
\end{picture}
\end{center}
\vspace{1.5in} \caption{The Cayley graph of $U_3(3)$ over
$PGL_2(14)$}
\end{figure}

The set of all double cosets [$w$]=${\cal N}w{\cal N}$, coset
stabilizing subgroups
${\cal N}^{(w)}$, and the number of single cosets each contains are exhibited in Table 1. \\ \\
\begin{table}
\begin{tabular}{l l r}
Label [$w$] & Coset stabilizing subgroup ${\cal N}^{(w)}$ & Number
\cr
 & & of cosets  \cr \hline
\\
$[\star]$ & Since $\cal N$ is transitive on T & 1 \cr
\\
$[{\bf 0}]$ & \parbox[t]{3 in}{\raggedright ${\cal N}^{({\bf 0})}
= N^{\bf 0}$ \mbox{$\langle y,({\bf 1},{\bf 4},{\bf 3},{\bf 5})
({\bf 2},{\bf 6})(1,4)(3,0,5,6)$},\mbox{$({\bf 1},{\bf 2},{\bf
3},{\bf 6})({\bf 4},{\bf 5})(1,2)(3,0,6,5),$} \mbox{$({\bf 2},{\bf
5})({\bf 4},{\bf 6})(2,4)(3,0)\rangle$} $\cong S_4$, has orbits
$\{{\bf 0}\}$,$\{1,2,14\}$, $\{3,5,6,0\}$ and $\{{\bf 1},{\bf
2},{\bf 3},{\bf 4},{\bf 5},{\bf 6}\}$ on T.} & 14 \cr
\\
$[{\bf 0}0]$=$[{\bf 0}]$
\\
& \parbox[t]{3 in}{\raggedright ${\bf 0}1 \sim 1{\bf 0} \sim {\bf
1}{\bf 3} \sim {\bf 3}{\bf 1} \sim 24 \sim 42 $} & \cr &
\parbox[t]{3 in}{\raggedright (Note
that $[{\bf 3}1] = [{\bf 0}1]$)} \\
$[{\bf 0}1] = [{\bf 0}{\bf 1}]$ & \parbox[t]{3 in}{\raggedright
$\forall i \in \{{\bf 1},{\bf 3},{\bf 0},{\bf 1},{\bf 2},{\bf
4}\}$} & 21 \cr
$[{\bf 0}1]=[7,i]$\\
& \parbox[t]{3 in}{\raggedright $\forall j \in$
\mbox{$\{2,4,5,6,10,12,13,14,15,16,17,18,19,20,$}
\mbox{$21,22,23,24,25,26,27,28,29,30,31,32,33,34,35,36\}$}} & \cr
$[7,8,j]=[7]$ \\ \hline
\end{tabular}
\caption{The double cosets [$w$]=${\cal N}w{\cal N}$ in $\cal G$.}
\end{table}
\\
\section{The algorithms and programs}
{\bf Algorithm I.} We start with the symmetric presentation of
\begin{equation} U_3(3):2\cong \frac{2^{*(7+7)}:PGL_2(7)}{t= {\bf  0} 0 {\bf
0}, y=({\bf 0} 1)^2}\end{equation}
in which the three generators x,
y and t correspond to the permutations a and c, and the symmetric
generator $t_{0}$, respectively. Next we compute the action of
$U_3(3):2$ on the 36 right cosets of subgroup $PGL_2(7)= <x,y>$ in
$U_3(3):2$, for example, the MAGMA command {\bf CosetAction}$(U3,
\mbox{sub}\langle  U3 | x, y \rangle )$ gives $J_{1}$ in its
action on 36 letters. We now form the symmetric generators ${\cal
T} = \{t_{i}; i=0,1,\ldots,14\}$ as permutations on 36 letters and
store them as $ts$, a sequence of length 14. Next we build $cst$,
a sequence of length 36 whose terms are sequences of integers,
representing words in the symmetric generators. These words form a
complete set of coset representatives for $PGL_2(7)$ in $U_3(3)$,
and correspond to the certain ordering, for example the ordering
determined by the MAGMA function {\bf CosetAction}. Represent
$U_3(3):2\cong \frac{2^{*(7+7)}:PGL_2(7)}{t= {\bf  0} 0 {\bf 0},
y=({\bf 0} 1)^2}$ as a permutaition group on the right cosets of
$PGL_2(7)$ in $U_3(3):2$ and consider a permutation $p$ of
$U_3(3):2$ on 36 letters. Convert $p$ to its canonical symmetric
representation form; that is, write it as a permutation of
$PGL_2(7)$ on 14
letters followed by a word of length of at most 2 in the symmetric generators.\\
{\bf Step I} $1^p =\cal{N}p = \cal{N}w$ $\implies$ $pw^-1 \in
\cal{N}$. $pw^-1$ is identified with an element of $\cal{N}$
by its action on the fourteen cosets whose representatives are of length 1.\\
{\bf Step II} $pw^-1$w is the permutation $p$ in its canonical
symmetric represntation form.\\
{\bf Algorithm II} Convert an element x of $U_3(3)$ given in the
symmetric representation form into a permutation on 36 letters.
Thus x is of the form $x = n w$, where $n \in \cal{N}$ and w a
word in the symmetric generators $t_i$s. Now the image of w is a
permutation on 36 letters and the image of $n$ is determined by
its action on the fourteen cosets whose representatives have
length 1.\\
{\bf Algorithm II} Given two elements of $U_3(3):2$ symmetrically
represented as $(xx,uu)$ and $(yy,vv)$, the procedure {\bf mult}
uses $ts$ and $cst$ to return the product $(zz,ww)$. As described
in the body of the paper the procedure {\bf cenelt} is used to
return generators for the centralizer of a given element
$(xx,uu)$, themselves symmetrically represented.
 As explained above, any element of $\cal G$ can
be written, not necessarily uniquely, as the product of a
permutation of ${\cal N} \cong PGL_{2}(14) $ followed by a word of
length at most four in the symmetric generators. In section $1.5$
of \cite{Cur2} the first author outlines a procedure for
multiplying elements represented in this fashion.It may be useful
to
see this process being carried out manually: \\
{\bf Example}
\begin{center}
\begin{tabular}{l c l}
$\pi_{{\bf 0}0}{\bf 12}\pi_{{\bf 2}3} {\bf 5}6$ & =
& $\pi_{{\bf 0}0}\pi_{{\bf 2}3}{\bf 12}^{\pi_{{\bf 2}3}}{\bf 5}6$ \\
& = & $\pi_{{\bf 0}0}\pi_{{\bf 2}3}03{\bf 5}6$  \\
& = &  $({\bf 1,4,3,2,6,5,0})(1,0,5,6,2,3,4)03{\bf 5}6$ \\
& = & $({\bf 1,4,3,2,6,5,0})(1,0,5,6,2,3,4)({\bf 2,4})({\bf 5,6})
(2,4)(5,6)30{\bf 5}6$\\
& = & $({\bf 1,2,5,0})({\bf 3,4})(1,0,6,4)(2,3)30{\bf 5}6$ \\
& = & $({\bf 1,2,5,0})({\bf 3,4})(1,0,6,4)(2,3)\pi_{{\bf 0}5}3^{\pi_{{\bf 0}5}}0{\bf 5}6$ \\
& = & $({\bf 1,2,5,0})({\bf 3,4})(1,0,6,4)(2,3)({\bf 1}, 3)({\bf
2}, 6)({\bf 3}, 2)
({\bf 4}, 4)({\bf 5}, 0)({\bf 6}, 1)({\bf 0}, 5){\bf 1}0{\bf 5}6$ \\
& = & $({\bf 1,2,5,0})({\bf 3,4})(1,0,6,4)(2,3)({\bf 1}, 3)({\bf
2}, 6)({\bf 3}, 2)({\bf 4}, 4)({\bf 5}, 0)({\bf 6}, 1)({\bf 0}, 5){\bf 1}0{\bf 5}6$ \\
& = & $({\bf 1,2,5,0})({\bf 3,4})(1,0,6,4)(2,3)({\bf 1}, 3)({\bf
2}, 6)({\bf 3}, 2)({\bf 4}, 4)({\bf 5}, 0)({\bf 6}, 1)({\bf 0}, 5){\bf 1}0{\bf 5}6$ \\
& = & $({\bf 1}, 6, {\bf 4}, 2)({\bf 2}, 0)({\bf 3}, 4, {\bf 6}, 1,
{\bf 5}, 5, {\bf 0}, 3){\bf 1}06$ \\
& = & $({\bf 1}, 6, {\bf 4}, 2)({\bf 2}, 0)({\bf 3}, 4, {\bf 6},1,{\bf 5}, 5, {\bf 0}, 3)
({\bf 3, 0})({\bf 5, 6})(2, 6)(4, 5)0{\bf 1}6$ \\
& = & $({\bf 1}, 2)({\bf 2}, 0)({\bf 3}, 5)({\bf 4}, 6)({\bf 5}, 4)({\bf 6}, 2)
({\bf 0}, 3)0{\bf 1}6$ \\
& = & $({\bf 1}, 2)({\bf 2}, 0)({\bf 3}, 5)({\bf 4}, 6)({\bf 5},
4)({\bf 6}, 2)
({\bf 0}, 3)\pi_{{\bf 1}6}0^{\pi_{{\bf 1}6}}{\bf 1}$ \\
& = & $({\bf 1}, 2)({\bf 2}, 0)({\bf 3}, 5)({\bf 4}, 6)({\bf 5},
4)({\bf 6}, 2)
({\bf 0}, 3)({\bf 1}, 6)({\bf 2}, 0)({\bf 3}, 1)({\bf 4}, 2)({\bf 5}, 3)({\bf 6}, 4)({\bf 0}, 5){\bf 2 1}$, \\
\end{tabular}
\end{center}
which is in canonically shortest form. In this paper we have
computerised such a procedure in two independent ways. The first,
see Appendix A, makes full use of the ease with which MAGMA
handles permutations of such low degree. Elements represented as
above are transformed into permutations on 266 letters, and any
group theoretic function can then be applied before transformation
back into the symmetric representation. Thus, for example, the
procedure {\em cenelt} returns generators for the centralizer of a
given element symmetrically represented. One can readily write
procedures to perform whatever task one chooses in MAGMA, but keep
a record of the results in this short form.

The programs given in Appendix B are rather more interesting, both
mathematically and computationally. To multiply two elements we
first use {\em unify} to express
\begin{displaymath}
{\pi}u.{\sigma}v = {\pi}{\sigma}.u^{{\sigma}}v
\end{displaymath}
where ${\pi},{\sigma} \in {\cal N}$ and $u,v$ are words in the
elements of ${\cal T}$, as a single sequence ss of length
$11+length(u)+length(v) \le 19$. The first 11 entries give the
permutation ${\pi}{\sigma}$ and the remainder represents a word of
length $\le$ 8 in the elements of ${\cal T}$. The procedure {\em
canon} now puts ss into its canonically shortest form. No other
representations of group elements are used; words in the symmetric
generators are simply shortened by application of the relations
(and their conjugates under ${\cal N}$).
Working interactively the response is immediate. \\
The first type of procedure is heavily MAGMA-dependent, but can be
readily modified for other packages such as GAP. On the other hand
the second type of procedure, although written in MAGMA here,
could have been written in any high-level
language. \\
This work represents a first step in a programme to provide
symmetric representations for interesting finite groups. Curtis
has obtained suitable symmetric presentations for several sporadic
simple groups; these will be computerised in a similar manner. Our
next aim is `nested' symmetric representations to deal with
much larger examples. \\  \\
\begin{center}
{\Large\bf Appendix}
\end{center}
\appendix
\section{Program 1}
\scriptsize We start with a presentation of $J_{1}$, based on the
symmetric presentation (1), in which the three generators x, y and
t correspond to the permutations a and c, and the symmetric
generator $t_{0}$, respectively. Thus the MAGMA command {\bf
CosetAction}$(J,  \mbox{sub}\langle  J | x, y \rangle )$ gives
$J_{1}$ in its action on 266 letters. We now form the symmetric
generators ${\cal T} = \{t_{i}; i=0,1,\ldots,X\}$ as permutations
on 266 letters and store them as $ts$, a sequence of length 11.
Next we build $cst$, a sequence of length 266 whose terms are
sequences of integers, representing words in the symmetric
generators. These words form a complete set of coset
representatives for $L_{2}(11)$ in $J_{1}$, and correspond to the
ordering determined by the MAGMA function {\bf CosetAction}. Given
two elements of $J_{1}$ symmetrically represented as $(xx,uu)$ and
$(yy,vv)$, the procedure {\bf mult} uses $ts$ and $cst$ to return
the product $(zz,ww)$. As described in the body of the paper the
procedure {\bf cenelt} is used to return generators for the
centralizer of a given element $(xx,uu)$, themselves symmetrically
represented.

\scriptsize \begin{verbatim}
U3fmt := recformat<
 /* Data structure for the symmetric representation of U_3(3). */
                   U3: GrpPerm,
                   PGL14: GrpPerm,
                   cst: SeqEnum,    // of integer sequences
                   ts: SeqEnum,     // of elements of U3
                   tra: SetIndx   // of elements of PGL14
>;

//-----------------------------------------------------------------------------
prodim := function(pt, Q, I)
 /*
 Return the image of pt under permutations Q[I] applied sequentially.
 */
                   v := pt;
                   for i in I do
                   v := v^(Q[i]);
                   end for;
                   return v;
 end function;

//-----------------------------------------------------------------------------
symrep := function()
 /*
 Initialize the data structures for the symmetric representation of U3.
 */
G<x,y,t,s>:=Group<x,y,t,s|x^7, y^2, t^2, (x^-1 * t)^2 , (y * x)^3
, t * x^-1 * y * x \
* t * y , x^2 * y * x^3 * y * x^-4 * y * x^-4 * y*x , s^2, (s^(x^3),y),(s^(x^4\
), x*y), t * s * s^t * s, y * (s * s^(t * x^6))^2>;

   // Construct the sequence of 14 symmetric generators ts as permutations
   // on 36 letters.

   f, U3, k := CosetAction(G, sub< G | x, y,t>);
   ts := [ Id(U3): i in [1 .. 14] ];
   for i in [1..7] do
   ts[i]:=(s^(x^i))@f;
   end for;
   for i in [8..14] do
   ts[i]:=((s^t)^(x^(14-i)))@f;
   end for;


   // Construct representatives cst for the control subgroup N = PGL(2, 14)
   // as words in the symmetric generators consisting of the empty word,
   // 14 words of length one and 21 words of length two.

   S14 := SymmetricGroup(14);
   aa:=S14!(1,2,3,4,5,6,7)(14,13,12,11,10,9,8);
   bb:=S14!(2,6)(4,5)(14,10)(13,12);
   cc:=S14!(7,14)(1,8)(2,9)(3,10)(4,11)(5,12)(6,13);
   PGL14 := sub< S14 | aa, bb, cc >;

   cst := [null : i in [1 .. 36]] where null is [Integers() | ];
   for i := 1 to 14 do
      cst[prodim(1, ts, [i])] := [i];
   end for;

   tra := Transversal(PGL14, sub<PGL14 | bb,cc^(aa*bb*aa^(-1)*bb)>);
   for i := 1 to 21 do
      ss := [7,1]^tra[i];
      cst[prodim(1, ts, ss)] := ss;
   end for;


   return rec<U3fmt |
      U3 := U3, PGL14 := PGL14, cst := cst, ts := ts, tra := tra>;
end function;

//-----------------------------------------------------------------------------


//-----------------------------------------------------------------------------
sym2per := function(U3Des, x)
  /*
Convert an element x of U3 in the symmetric repesentation into a
permutation acting on 36 letters. The image of an element of N is
determined by its action on the fourteen cosets whose
representatives have length one.
 */
           U3 := U3Des`U3; cst := U3Des`cst; ts := U3Des`ts;  PGL14 := U3Des`PGL14;
           tra := U3Des`tra;

           xx := U3Des`PGL14 ! x[1]; uu := x[2];
           p := [1 : i in [1 .. 36]];
           for i := 1 to 14 do
           p[prodim(1, ts, [i])] := prodim(1, ts, [i]^xx);
           end for;
           for i := 1 to 21 do
           t := [7,1]^tra[i];
           p[prodim(1, ts, t)] := prodim(1, ts, t^xx);
           end for;
           return (U3 ! p) * &*[U3|ts[uu[j]]: j in [1 .. #uu]];
end function;

//-----------------------------------------------------------------------------
per2sym := function(U3Des, p) /* Convert permutation p of U3 on 36
letters into its symmetric representation. The image of 1 under p
gives the coset representative for Np as a word ww in the
symmetric generators. Multiplication of p by the symmetric
generators of ww in reverse order yields a permutation which can
be identified with an element of N by its action on the 14 cosets
of length one. */
           U3 := U3Des`U3; cst := U3Des`cst; ts := U3Des`ts; PGL14 := U3Des`PGL14;
           ww := cst[1^p];
           tt := p * &*[U3|ts[ww[#ww - l + 1]]: l in [1 .. #ww]];
           zz := PGL14![rep{j: j in [1..14] | (1^ts[i])^tt eq 1^ts[j]}: i in [1..14]];
           return <zz, ww>;
end function;

//-----------------------------------------------------------------------------
cenelt := function(U3Des, x)
/*
 Construct the centralizer of element x of J1 given in its symmetric representation.
 An example of how all the standard procedures of MAGMA can be utilized by:
 transformation to permutations, application of the procedure,
 transformation back to symmetric representation.
*/
           cent := Centralizer(U3Des`U3, sym2per(U3Des, x));
           return <Order(cent), [per2sym(U3Des, c): c in Generators(cent)]>;
end function;
\end{verbatim}
\newpage
\section{Program 2}
In this program we assume detailed knowledge of the {\em control
subroup} \mbox{${\cal N} \cong PGL_{2}(7)$}, but use no
representation of elements of $J_{1}$ other than their symmetric
representation. Firstly, the procedure {\em unify} uses the
identity ${\pi}u.{\sigma}v = {\pi}{\sigma}.u^{{\sigma}}v$ to
 combine two symmetrically represented
elements $(xx,uu)$ and $(yy,vv)$ into a single sequence $ss$ of
length $(11+length(uu)+length(vv))$, which represents a
permutation of ${\cal N}$ followed by a word of length $\le 8$ in
the elements of ${\cal T}$. The procedure {\em canon} then takes
such a sequence
and reduces it to its shortest form using the following recursive algorithm. We make use of the relations \mbox{$t_it_jt_it_jt_i = \pi_{ij}$} \  and  \mbox{$t_it_jt_kt_it_j = \sigma_{ijk}$}, for $\{ i,j,k \}$ a {\em special triple} i.e. a triple fixed pointwise by an involution of $\cal N$.\\
{\em The algorithm}: \\
{\bf Step I:} If two adjacent symmetric generators are equal,
delete them.\\
{\bf Step II:} If a string $t_{i}t_{j}t_{i}$ appears, replace it
by ${\pi}_{ij}t_{i}t_{j}$ and move the permutation ${\pi}_{ij}$
over the preceding symmetric generators in the standard manner. \\
{\bf Step III:} If a string $t_{i}t_{j}t_{k}$ appears with
$\{i,j,k\}$ a special triple replace it by
${\sigma}_{ijk}t_{j}t_{i}$,
and move the permutation to the left as above.\\
Having completed the above, if $length(ss) \le 11+3=14$,
finish. Otherwise we may assume all strings $t_it_jt_k$ have $\{i,j,k \}$ non-special. \\
{\bf Step IV:}. For each string $t_it_jt_kt_l$ construct a permutation $\sigma$ and a 5-cycle $\rho$ such that $\sigma t_it_jt_k$ commutes with $\rho$. [In fact $\sigma = \sigma_{pqr}$ where $p=i^{\pi_{jk}}, q=k^{\pi_{i,j}},r=j^{\pi_{ki}}$]. If $l$ is not fixed by $\rho$ we replace $t_it_jt_k$ by $\sigma^{-1}(\sigma t_it_jt_k)^{\rho^m}$, where $m$ is such that $l \in \{ i,j,k \}^{\rho^m}$. If $length(ss) \le 15$, finish. Otherwise: \\
{\bf Step V:} We may now assume that any string of length 5 is one
of the eleven special pentads preserved by $L_2(11)$. We use the
identity $267X8 = \sigma_{X64} X26$ from \cite[page 304]{Cur2}.
i.e. for [$i,j,k,l,m$] an even permutation of a special pentad
$$
t_it_jt_kt_lt_m = \sigma_{ljn} t_lt_it_j,
$$
where $n = l^{\pi_{ij}}$. \\
After each {\bf Step} recall {\em canon}. \\
\par
\scriptsize \begin{verbatim}
 /*
 Define the projective general
 linear group PGL(2,14) as permutations of degree 14
 */
S14:=SymmetricGroup(14);
a:=S14!(1,2,3,4,5,6,7)(14,13,12,11,10,9,8);
b:=S14!(2,6)(4,5)(14,10)(13,12);
c:=S14!(7,14)(1,8)(2,9)(3,10)(4,11)(5,12)(6,13);
N:=sub<S14|a,b,c>;
 /*--------------------------------------------------------------
Coset Stabiliser N(0,1) /* NS71:=sub<N|N!(2, 6)(4, 5)(10, 14)(12,
13), N!(1, 13, 7, 12)(2, 9, 6, 10, 5, 11, 4, 14)(3, 8)>;
NS71:=sub<N|b,c^(a*b*a^(-1)*b)>;
 /*-------------------------------------------------------------
 Six names for the 21-orbit
 /*
ss:=[[7,1],[1,7],[13,12],[12,13],[3,8],[8,3]];
sp:=[Id(N),N!(2,5)(4,6)(14,10)(9,11),N!(7,1)(2,6,5,4)(14,9,10,11)(12,13),
N!(7,1)(2,4,5,6)(14,11,10,9)(12,13),N!(2,6)(4,5)(14,10)(12,13),
N!(2,4)(5,6)(9,11)(12,13)];
trans1:=Transversal(N,Stabiliser(N,{7,14})); prs:={@ {7,14}^x : x
in trans1 @}; pi:=N!(7,14)(1,8)(2,9)(3,10)(4,11)(5,12)(6,13);
pis:=[ pi^x : x in trans1]; trans2:=Transversal(N,NS71);

 /*
----------------------------------------------------------------------------
 Given two symmetrically represented elements of J1, where x and y
 are permutations of N and u and v are words in the symmetric
 generators, return a single sequence of length 11 + l(u) + l(v)
 using the above identity.
 */
 Unify := func< x, y, u, v | [ y[x[i]] : i in [1..#x] ] cat y[15..#y]
                             where x is Eltseq(x) cat u
                             where y is Eltseq(y) cat v >;
aaa:=N!(2,6)(4,5)(14,10)(13,12);
bbb:=N!(7,14)(1,8)(2,9)(3,10)(4,11)(5,12)(6,13);
Unify(aaa,bbb,u,v);
sss:=[8,9,10,11,12,13,14,1,2,3,4,5,6,7,7,14,7,1,8];
/*----------------------------------------------------------------------------
Return permutation of N given by the word ijiji in the symmetric
generators. */
Pi := func< i, j | pis[Index(prs,{i,j})]>;
 /*----------------------------------------------------------------------------
For ss a sequence representing a permutation of N followed by a
word in the symmetric generators, return an equivalent sequence of
canonically shortest length.
*/
function Canon(sss)
    s:=sss;
   // Step I.
  if exists(i){ i : i in [15..#s-1] | s[i] eq s[i+1] } then
        s := Canon( s[1..i-1] cat s[i+2..#s] );
  end if;
  // Step II.
  if exists(i){ i : i in [15..#s-1] | Index(prs,{s[i],s[i+1]}) ne 0} then
  s := Canon( [p[s[k]] : k in [1..i-1] ] cat s[i..i] cat s[i+2..#s]
              where p is Eltseq(Pi(s[i], s[i+1])) );
  end if;
  // Step III.
  if #s ge 17 then for x in trans2 do if exists(i){i : i in [1..6] |
  [s[15],s[16]]^x eq ss[i]} then if exists(j){j : j in [1..6] |
  (s[17])^x eq ss[j][2]} then s:=Canon([p[s[k]]  : k in [1..14]] cat
  (ss[j])^(x^(-1)) cat s[17..#s] where p is
  Eltseq((sp[i]^(-1)*sp[j])^(x^(-1)))); end if; end if; end for;end
  if;
  return s;
end function;

S14:=SymmetricGroup(14);
a:=S14!(1,2,3,4,5,6,7)(14,13,12,11,10,9,8);
b:=S14!(2,6)(4,5)(14,10)(13,12);
c:=S14!(7,14)(1,8)(2,9)(3,10)(4,11)(5,12)(6,13);
N:=sub<S14|a,b,c>; NS71:=sub<N|N!(2, 6)(4, 5)(10, 14)(12, 13),
N!(1, 13, 7, 12)(2, 9, 6, 10, 5,\
 11, 4, 14)(3, 8)>;
ss:=[[7,1],[1,7],[13,12],[12,13],[3,8],[8,3]];
sp:=[Id(N),N!(2,5)(4,6)(14,10)(9,11),N!(7,1)(2,6,5,4)(14,9,10,11)(12,13),
N!(7,1)(2,4,5,6)(14,11,10,9)(12,13),N!(2,6)(4,5)(14,10)(12,13),
N!(2,4)(5,6)(9,11)(12,13)];
trans1:=Transversal(N,Stabiliser(N,{7,14})); prs:={@ {7,14}^x : x
in trans1 @}; pi:=N!(7,14)(1,8)(2,9)(3,10)(4,11)(5,12)(6,13);
pis:=[ pi^x : x in trans1]; Pi := func< i, j |
pis[Index(prs,{i,j})]>; trans2:=Transversal(N,NS71); function
Canon(sss)
    s:=sss;
   // Step I.
  if exists(i){ i : i in [15..#s-1] | s[i] eq s[i+1] } then
        s := Canon( s[1..i-1] cat s[i+2..#s] );
   end if;
    // Step II.
 if exists(i){ i : i in [15..#s-1] | Index(prs,{s[i],s[i+1]}) ne 0} then
 s := Canon( [ p[s[k]] : k in [1..i-1] ] cat s[i..i] cat s[i+2..#s]
                  where p is Eltseq(Pi(s[i], s[i+1])) );
end if;
// Step III.
if #s ge 17 then if exists(i,j,x){<i,j,x> : i in [1..6], j in
[1..6] ,  x in trans2 | ss[i]^x eq [s[15],s[16]] and ss[j][2]^x eq
s[17]} then s:=Canon([p[s[k]]  : k in [1..14]] cat [ss[j][1]^x]
cat s[18..#s] where p is Eltseq((sp[i]^(-1)*sp[j])^x)); else if
exists(i,j,x){<i,j,x> : i in [1..6], j in [1..6] ,  x in trans2 |
ss[i]^x eq [s[15],s[16]] and Index(prs,{ss[j][2]^x,s[17]}) ne 0}
then s:=Canon([p[s[k]]  : k in [1..14]] cat
[ss[j][1]^(Pi(ss[j][2]^x,s[17])) ] cat [ss[j][2]^x] cat s[18..#s]
where p is Eltseq((sp[i]^(-1)*sp[j])^x*Pi(ss[j][2]^x,s[17]))); end
if; end if; end if;
return s;
end function;
 /*------------------------------------------------------------------------------
 Return the product of two symmetrically represented elements of J1
 */ Prod := function( x, u, y, v)
           t := Canon(Unify(x, u, y, v));
           return N!t[1..14], t[15..#t];
    end function;

/*--------------------------------------------------------------------------
 Return the inverse of a symmetrically represented element of J1
*/
   Invert := func< x, u | x^-1, [u[#u-i+1]^(x^-1) : i in [1..#u]] >;


\end{verbatim}

\end{document}